\documentclass[a4paper, 10pt]{amsart}
\usepackage{amssymb}
\usepackage[utf8]{inputenc}
\usepackage[french, english]{babel}
\usepackage[T1]{fontenc}
\usepackage[margin=2cm]{geometry}
\usepackage{comment}
\usepackage{hyperref}

\title{A note on weak amenability for free products of discrete quantum groups}
\author{Amaury Freslon}
\email{freslon@math.jussieu.fr}
\urladdr{http://www.math.jussieu.fr/~freslon}
\date{}
\address{Univ. Paris Diderot, Paris Cit\'e Sorbonne, IMJ, UMR 7586, 175 rue du Chevaleret, 75013, Paris, France}

\theoremstyle{plain}
\newtheorem{theorem}{Theorem}[section]
\newtheorem{e-proposition}[theorem]{Proposition}

\newtheorem{lemma}[theorem]{Lemma}

\theoremstyle{definition}
\newtheorem{e-definition}[theorem]{Definition}

\theoremstyle{remark}

\newtheorem{example}{Example}

\newcommand{\B}{\mathcal{B}}
\newcommand{\G}{\mathbb{G}}
\newcommand{\h}{\widehat}
\renewcommand{\i}{\imath}
\newcommand{\w}{\omega}

\begin{document}

\begin{abstract}
We prove that the Cowling-Haagerup constant of a reduced free product of weakly amenable discrete quantum groups with Cowling-Haagerup constant equal to 1 is again equal to 1.
\end{abstract}

\maketitle

\section{Introduction}

In geometric group theory, weak amenability is an approximation property which is satisfied by a large class of groups (for example free or even hyperbolic groups \cite{ozawa2007weak}) but is strong enough to give interesting properties for the von Neumann algebras associated to these groups (for example related to deformation/rigidity techniques \cite{ozawa2010class}, \cite{ozawa2010classbis}). The stability of this property under free products is still an open problem. However, using a very general version of the Khintchine inequality, E. Ricard and Q. Xu were able to prove in \cite{ricard2005khintchine} that if $(G_{i})$ is a family of weakly amenable discrete groups \emph{with Cowling-Haagerup constant equal to $1$}, then their free product is again weakly amenable, and its Cowling-Haagerup constant is also equal to $1$. The proof uses a classical characterization of the bounded functions on a group giving rise to completely bounded multipliers. This characterization having recently been generalized to arbitrary locally compact quantum groups by M. Daws in \cite{daws2011multipliers}, we are able to prove an analogue of Ricard and Xu's result in the setting of discrete quantum groups.

We are very grateful to P. Fima for many useful discussions and advice and to E. Blanchard and R. Vergnioux for their careful proof-reading of an early version. We also thank M. Brannan for pointing to us the paper \cite{daws2011multipliers}. 

\section{Definitions and notations}

In this section, we introduce the notations and results about discrete quantum groups which will be used in this paper. For two Hilbert spaces $H$ and $K$, $\mathcal{B}(H, K)$ will denote the set of bounded linear maps from $H$ to $K$ and $\mathcal{B}(H)=\mathcal{B}(H, H)$. In the same way, we use the notations $\mathcal{K}(H, K)$ and $\mathcal{K}(H)$ for compact linear maps. We will denote by $\mathcal{B}(H)_{*}$ the predual of $\mathcal{B}(H)$, i.e. the Banach space of all normal linear forms on $\mathcal{B}(H)$. On any tensor product $A\otimes B$, we define the flip operator $\Sigma : x\otimes y \mapsto y\otimes x$. We will use the usual leg-numbering notations : for an operator $X$ acting on a tensor product we set $X_{12}:=X\otimes1$, $X_{23}:=1\otimes X$ and $X_{13}:=(\Sigma\otimes 1)(1\otimes X)(\Sigma\otimes 1)$. For a subset $B$ of a topological vector space $C$, $\overline{\rm{span}} B$ will denote the \emph{closed linear span} of $B$ in $C$. The symbol $\otimes$ will denote the \emph{minimal} (or spatial) tensor product of C*-algebras or the topological tensor product of Hilbert spaces. Let $A$ be a C*-algebra together with a distinguished state $\varphi$. Taking the completion of $A$ with respect to the scalar product $\langle a, b\rangle = \varphi(a^{*}b)$ yields a Hilbert space which will be denoted $L^{2}(A, \varphi)$. If the scalar product is $\langle a, b\rangle^{op} = \varphi(ab^{*})$, then the completion will be denoted $L^{2}(A, \varphi)^{op}$.

Let $\mathbb{G} = (C(\mathbb{G}), \Delta))$ be a compact quantum group in the sense of \cite{woronowicz1995compact} and let $\widehat{\mathbb{G}} = (C_{0}(\widehat{\mathbb{G}}), \widehat{\Delta})$ be its dual discrete quantum group. The Haar state of $\mathbb{G}$ will be denoted $h$. We will always assume it to be faithful and identify $C(\mathbb{G})$ with its image under the GNS map. We will denote by $L^{2}(\mathbb{G})$ the Hilbert space of the GNS construction, by $\xi_{h}$ the cyclic separating vector and by $W$ the unitary operator on $L^{2}(\mathbb{G})\otimes L^{2}(\mathbb{G})$ defined by $W^{*}(\xi\otimes a\xi_{h}) = \Delta(a)(\xi\otimes \xi_{h})$ for $\xi \in L^{2}(\mathbb{G})$ and $a\in C(\mathbb{G})$. Then, $W$ is a \emph{multiplicative unitary} in the sense of \cite{baaj1993unitaires}, i.e. $W_{12}W_{13}W_{23} = W_{23}W_{12}$. We have the following equalities :
\begin{eqnarray*}
C(\mathbb{G}) & = & \overline{\rm{span}}(\i\otimes \mathcal{B}(L^{2}(\mathbb{G}))_{*})(W) \\
\Delta(x) & = & W^{*}(1\otimes x)W \\
C_{0}(\widehat{\mathbb{G}}) & = & \overline{\rm{span}}(\mathcal{B}(L^{2}(\mathbb{G}))_{*}\otimes id)(W) \\
\widehat{\Delta}(x) & = & \Sigma W(x\otimes 1)W^{*}\Sigma.
\end{eqnarray*}
These faithful representations on $L^{2}(\mathbb{G})$ allow us to define the von Neumann algebras $L^{\infty}(\mathbb{G}) = C(\mathbb{G})''$ and $\ell^{\infty}(\widehat{\mathbb{G}})=C_{0}(\widehat{\mathbb{G}})''$, and one can check that $W\in L^{\infty}(\mathbb{G})\otimes\ell^{\infty}(\widehat{\mathbb{G}})$. Finally, we set $\widehat{W} = \Sigma W^{*}\Sigma$.

\section{Preliminaries}

Let $\widehat{\mathbb{G}}$ be a discrete quantum group and $a\in \ell^{\infty}(\widehat{\mathbb{G}})$. The \emph{left multiplier} associated to $a$ is the map $m_{a}$ on $(\i\otimes \mathcal{B}(L^{2}(\mathbb{G}))_{*})(W)$ defined by $(m_{a}\otimes id)(W) = (1\otimes a)W$. A net $(a_{i})$ of elements of $\ell^{\infty}(\widehat{\mathbb{G}})$ is said to \emph{converge pointwise} to $a\in \ell^{\infty}(\widehat{\mathbb{G}})$ if $a_{i}p\rightarrow ap$ in $\mathcal{B}(L^{2}(\mathbb{G}))$ for any minimal central projection $p\in \ell^{\infty}(\widehat{\mathbb{G}})$. An element $a\in \ell^{\infty}(\widehat{\mathbb{G}})$ is said to have \emph{finite support} if $ap$ vanishes for all but finitely many central projections $p\in \ell^{\infty}(\widehat{\mathbb{G}})$.

For a discrete group $G$, it is known that a bounded function $\varphi : G\rightarrow \mathbb{C}$ gives rise to a completely bounded multiplier if and only if there exists a Hilbert space $K$ and two families $(\xi_{s})_{s\in G}$ and $(\eta_{t})_{t\in G}$ of vectors in $K$ such that $\varphi(s) = \langle \eta_{t}, \xi_{st}\rangle$ (which is usually written $\varphi(st^{-1}) = \langle \eta_{t}, \xi_{s}\rangle$). Moreover, the completely bounded norm is then $\inf \|(\xi_{s})\|_{\infty}\|(\eta_{t})\|_{\infty}$ (see e.g. \cite{bozejko1984herz}). The following theorem \cite[Prop 4.1 and Thm 4.2]{daws2011multipliers} gives the quantum analogue of this characterization.

\begin{theorem}[Daws \cite{daws2011multipliers}]\label{thm:quantumgilbert}
Let $\widehat{\mathbb{G}}$ be a discrete quantum group and $a\in \ell^{\infty}(\widehat{\mathbb{G}})$. Then $m_{a}$ extends to a competely bounded multiplier on $\mathcal{B}(L^{2}(\mathbb{G}))$ if and only if there exists a Hilbert space $K$ and two maps $\alpha, \beta \in \mathcal{B}(L^{2}(\mathbb{G}), L^{2}(\mathbb{G})\otimes K)$ such that $(1\otimes \beta)^{*}\widehat{W}_{12}^{*}(1\otimes \alpha)\widehat{W} = a\otimes 1.$ Moreover, we then have $m_{a}(x) = \beta^{*}(x\otimes 1)\alpha$ and we can chose $\alpha$ and $\beta$ to have norm equal to $\sqrt{\|m_{a}\|_{cb}}$.
\end{theorem}

\begin{proof}
We only give the construction of $\alpha$ and $\beta$ since we will need their precise form later on. Assume $m_{a}$ to be completely contractive. By Wittstock's factorization theorem, there is a representation $\pi : \mathcal{B}(L^{2}(\mathbb{G})) \rightarrow \mathcal{B}(K)$ and two isometries $P, Q \in \mathcal{B}(L^{2}(\mathbb{G}), K)$ such that for all $x\in \mathcal{B}(L^{2}(\mathbb{G}))$, $m_{a}(x) = Q^{*}\pi(x)P$. Then, setting $\alpha = (\i\otimes \pi)(\widehat{W})(1\otimes P)\widehat{W}(\i\otimes \xi_{h})$ and $\beta = (\i\otimes \pi)(\widehat{W})(1\otimes Q)\widehat{W}(\i\otimes \xi_{h})$ yields the result.
\end{proof}

From this we easily deduce that a multiplier $m_{a}$ on a discrete quantum group $\widehat{\mathbb{G}}$ has a completely bounded extension to $\mathcal{B}(L^{2}(\mathbb{G}))$ if and only if it has a completely bounded extension to $C(\mathbb{G})$ or to $L^{\infty}(\mathbb{G})$ and that the completely bounded norms of these extensions are all equal. This justifies the following definition of weak amenability for discrete quantum groups.

\begin{e-definition}\label{de:quantumwa}
A discrete quantum group $\widehat{\mathbb{G}}$ is said to be \emph{weakly amenable} if there exists a net $(a_{i})$ of elements of $\ell^{\infty}(\widehat{\mathbb{G}})$ such that $a_{i}$ has finite support for all $i$, the net $(a_{i})$ converges pointwise to $1$ and the maps $m_{a_{i}}$ satisfy $\limsup \|m_{a_{i}}\|_{cb} < \infty$. The infimum of $\limsup \|m_{a_{i}}\|_{cb, \mathcal{B}(L^{2}(\mathbb{G}))}$ for all nets satisfying the above conditions is denoted $\Lambda_{cb}(\widehat{\mathbb{G}})$ and called the \emph{Cowling-Haagerup constant} of $\widehat{\mathbb{G}}$. By convention, $\Lambda_{cb}(\widehat{\mathbb{G}})=\infty$ if $\widehat{\mathbb{G}}$ is not weakly amenable.
\end{e-definition}

\section{Free products of weakly amenable discrete quantum groups}

Recall that given two discrete quantum groups $\widehat{\mathbb{G}}$ and $\widehat{\mathbb{H}}$, there is a unique compact quantum group structure on the reduced free product $C(\mathbb{G})\ast C(\mathbb{H})$ with respect to the Haar states which is compatible with the canonical inclusions (it is defined in \cite{wang1995free}). By analogy with the classical case, the dual of this compact quantum group will be called the \emph{reduced free product} of $\widehat{\mathbb{G}}$ and $\widehat{\mathbb{H}}$ and denoted $\widehat{\mathbb{G}}\ast\widehat{\mathbb{H}}$. We will now prove that a reduced free product of weakly amenable discrete quantum groups with Cowling-Haagerup constant equal to $1$ has Cowling-Haagerup constant equal to $1$. This result has been proved in the classical case by E. Ricard and Q. Xu \cite[Theorem 4.3]{ricard2005khintchine} using the following result \cite[Prop 4.11]{ricard2005khintchine}.

\begin{theorem}[Ricard, Xu \cite{ricard2005khintchine}]\label{thm:freecmap}
Let $(B_{i}, \psi_{i})_{i\in I}$ be unital C*-algebras with distinguished states $(\psi_{i})$ having faithful GNS constructions. Let $A_{i}\subset B_{i}$ be unital C*-subalgebras such that the states $\varphi_{i} =\psi_{i\vert A_{i}}$ also have faithful GNS construction. Assume that for each $i$, there is a net of finite rank maps $(V_{i, j})$ on $A_{i}$ converging to the identity pointwise, preserving the state and such that $\limsup_{j}\|V_{i, j}\|_{cb}=1$. Assume moreover that for each pair $(i, j)$, there is a completely positive unital map $U_{i, j} : A_{i} \rightarrow B_{i}$ preserving the state and such that $\| V_{i, j} - U_{i, j}\|_{cb} + \| V_{i, j}-U_{i, j}\|_{\mathcal{B}(L^{2}(A_{i}, \varphi_{i}), L^{2}(B_{i}, \psi_{i}))} + \| V_{i, j} - U_{i, j}\|_{\mathcal{B}(L^{2}(A_{i}, \varphi_{i})^{op}, L^{2}(B_{i}, \psi_{i})^{op})} \rightarrow 0.$ Then, the reduced free product of the family $(A_{i}, \varphi_{i})$ has Cowling-Haagerup constant equal to $1$.
\end{theorem}

\begin{theorem}\label{thm:quantumfreeproduct}
Let $(\widehat{\mathbb{G}}_{i})_{i\in I}$ be a family of discrete quantum groups with Cowling-Haagerup constant equal to $1$. Then $\Lambda_{cb}(\ast_{i}\widehat{\mathbb{G}}_{i})=1$.
\end{theorem}

The proof of this theorem is quite involved. In order to make it more clear, we will divide it into several lemmata, most of which are rather technical. We first introduce some general notations. Let $\h{\mathbb{G}}$ be a discrete quantum group, let $0 < \epsilon < 1$ and let $a\in \ell^{\infty}(\h{\mathbb{G}})$ be such that $\|m_{a}\|_{cb}\leqslant 1+\epsilon$. Let $\xi, \eta : L^{2}(\G)\rightarrow L^{2}(\G)\otimes K$ denote the two maps given by Theorem \ref{thm:quantumgilbert} with $\|\xi\| = \|\eta\|\leqslant \sqrt{1+\epsilon}$ and set $\gamma = (\xi + \eta)/2$ and $\delta = (\xi-\eta)/2$. Obtaining an approximation by a unital completely positive map is the first step.

\begin{lemma}\label{lem:cbapproximation}
Assume that $a = \widehat{S}(a)^{*}$ and that $m_{a}(1) = 1$. Then, there exists a u.c.p. map on $\mathcal{B}(L^{2}(\mathbb{G}))$ approximating $m_{a}$ up to $6\epsilon$ in completely bounded norm.
\end{lemma}

\begin{proof}
We know from \cite[Prop 2.6]{kraus1999approximation} that for any $x\in C(\mathbb{G})$, $\xi^{*}(x\otimes 1)\eta = m_{a}(x^{*})^{*} = m_{\h{S}(a)^{*}}(x)$. Thus,
\begin{equation*}
m_{a}(x) = \frac{1}{2}(m_{a}(x) + m_{\h{S}(a)^{*}}(x)) = \frac{1}{2}((\eta^{*}(x\otimes 1)\xi + \xi^{*}(x\otimes 1)\eta)) = M_{\gamma}(x) - M_{\delta}(x),
\end{equation*}
where $M_{\gamma}(x) = \gamma^{*}(x\otimes 1)\gamma$ and $M_{\delta}(x) = \delta^{*}(x\otimes 1)\delta$. The maps $M_{\gamma}$ and $M_{\delta}$ are completely positive, thus $\|M_{\gamma}\|_{cb} = \|\gamma\|^{2}\leqslant 1+\epsilon$ and evaluating at $1$ gives $\|1 + \delta^{*}\delta\| \leqslant 1+\epsilon$, i.e. $\|M_{\delta}\|_{cb} = \|\delta^{*}\delta\|\leqslant \epsilon$. We now want to perturb $M_{\gamma}$ into a \emph{unital} completely positive map. To do this, first note that $\|1 - \gamma^{*}\gamma\| = \|\delta^{*}\delta\|\leqslant \epsilon <1$, which implies that $\gamma^{*}\gamma$ is invertible, and set $\tilde{\gamma} = \gamma\vert \gamma\vert^{-1}$ where $\vert \gamma \vert = (\gamma^{*}\gamma)^{1/2}$. Note that $\|\tilde{\gamma} - \gamma\|\leqslant \epsilon$.
Thus, $M_{\tilde{\gamma}}$ is a unital completely positive map and
\begin{eqnarray*}
\|M_{\tilde{\gamma}} - M_{\gamma}\|_{cb} & = & \|M_{\gamma + (\tilde{\gamma}- \gamma)} - M_{\gamma}\|_{cb} \\
& \leqslant & \|\tilde{\gamma}-\gamma\|\|\gamma\| + \|\tilde{\gamma}-\gamma\|\|\gamma\| + \|\tilde{\gamma}-\gamma\|\|\tilde{\gamma}-\gamma\| \\
& \leqslant & \epsilon(2+3\epsilon)\leqslant 5\epsilon.
\end{eqnarray*}
This proves that $M_{\tilde{\gamma}}$ is a unital completely positive map approximating $m_{a}$ on $C(\mathbb{G})$ up to $6\epsilon$ in completely bounded norm.
\end{proof}

Set $D = \mathcal{B}(L^{2}(\mathbb{G}))$. We now want to prove that the previous approximation also works when the maps are seen as operators on $L^{2}(D, \tau)$ and $L^{2}(D, \tau)^{op}$, where $\tau(x) = \langle \xi_{h}, x(\xi_{h})\rangle$. Let us start with a purely computational lemma.

\begin{lemma}\label{lem:computational}
For any $\zeta\in K$, $(\i\otimes \pi)(\h{W})(\xi_{h}\otimes \zeta) = \xi_{h}\otimes \zeta$.
\end{lemma}

\begin{proof}
Note that since by definition $W(\xi\otimes \xi_{h}) = \xi\otimes \xi_{h}$ for any $\xi\in H$, we also have $\h{W}(\xi_{h}\otimes \xi) = \xi_{h}\otimes \xi$. For any $\theta_{1}, \theta_{2}, \xi\in L^{2}(\mathbb{G})$,
\begin{eqnarray*}
\langle \xi, (\i \otimes \w_{\theta_{1}, \theta_{2}})(\h{W})\xi_{h}\rangle & = & \langle  \xi\otimes \theta_{1}, \h{W}(\xi_{h}\otimes\theta_{2})\rangle \\
& = & \langle \xi, \xi_{h}\rangle \langle \theta_{1}, \theta_{2}\rangle \\
& = & \w_{\theta_{1}, \theta_{2}}(1)\langle \xi, \xi_{h}\rangle.
\end{eqnarray*}
Thus by density, we have $\langle \xi, (\i\otimes \w)(\h{W})\xi_{h}\rangle = \w(1)\langle \xi, \xi_{h}\rangle$ for any $\w\in \mathcal{B}(L^{2}(\mathbb{G}))_{*}$. Now, let $\zeta_{1}, \zeta_{2}\in K$ and $\xi\in L^{2}(\mathbb{G})$, then
\begin{eqnarray*}
\langle \xi\otimes \zeta_{1}, (\i\otimes \pi)(\h{W})(\xi_{h}\otimes \zeta_{2})\rangle & = & \langle \xi, (\i\otimes \w_{\zeta_{1}, \zeta_{2}}\circ\pi)(\h{W})\xi_{h}\rangle \\
& = & \w_{\zeta_{1}, \zeta_{2}}(\pi(1))\langle \xi, \xi_{h}\rangle \\
& = & \langle \xi\otimes \zeta_{1}, \xi_{h}\otimes \zeta_{2}\rangle.
\end{eqnarray*}
\end{proof}

This gives us a systematic way to investigate the $L^{2}$-norm of some specific operators.

\begin{lemma}\label{lem:technical1}
Let $T$ be any bounded linear operator from $L^{2}(\mathbb{G})$ to $K$ and set
\begin{equation*}
A(T) = (\i\otimes\pi)(\widehat{W})^{*}(1\otimes T)\widehat{W}(\i \otimes \xi_{h})\in \mathcal{B}(L^{2}(\mathbb{G}), L^{2}(\mathbb{G})\otimes K).
\end{equation*}
and $M_{A(T)}(x) = A(T)^{*}(x\otimes 1)A(T)$. Then, $\tau(M_{A(T)}(x^{*}x))\leqslant \|T\|^{2}\tau(x^{*}x)$ and $M_{A(T)}$ is a bounded operator on $L^{2}(D, \tau)$ of norm less than $\|T\|^{2}$. If moreover $A(T)^{*}A(T)$ is invertible, then $M_{A(T)\vert A(T)\vert^{-1}}$ is $\tau$-invariant.
\end{lemma}

\begin{proof}
Let us compute
\begin{eqnarray*}
A(T)\xi_{h} & = & (\i\otimes\pi)(\h{W})^{*}(1\otimes T)\h{W}(\xi_{h}\otimes \xi_{h}) \\
& = & (\i\otimes\pi)(\h{W})^{*}(1\otimes T)(\xi_{h}\otimes \xi_{h}) \\
& = & (\i\otimes\pi)(\h{W})^{*}(\xi_{h}\otimes T(\xi_{h})) \\
& = & \xi_{h}\otimes T(\xi_{h}).
\end{eqnarray*}
From this, we get
\begin{equation*}
\langle \xi_{h}, A(T)^{*}(x\otimes 1)A(T)\xi_{h}\rangle = \langle A(T)\xi_{h}, (x\otimes 1)A(T)(\xi_{h}) \rangle = \langle \xi_{h}, x(\xi_{h})\rangle \|T(\xi_{h})\|^{2}
\end{equation*}
and Kadison's inequality yields
\begin{eqnarray*}
\tau(M_{A(T)}(x)^{*}M_{A(T)}(x)) & \leqslant & \|A(T)\|^{2} \tau(M_{A(T)}(x^{*}x)) \\
& \leqslant & \|T\|^{2}\|A(T)\|^{2} \tau(x^{*}x) \\
& \leqslant & \|T\|^{4} \tau(x^{*}x).
\end{eqnarray*}

Let us now turn to $A(T)^{*}A(T)$. First,
\begin{eqnarray*}
A(T)^{*}A(T)\xi_{h} & = & (\i\otimes \xi_{h}^{*})\h{W}^{*}(1\otimes T^{*})(\i\otimes \pi)(\h{W})(\xi_{h}\otimes T(\xi_{h})) \\
& = & (\i\otimes \xi_{h}^{*})\h{W}^{*}(1\otimes T^{*})(\xi_{h}\otimes T(\xi_{h})) \\
& = & (\i\otimes \xi_{h}^{*})\h{W}^{*}(\xi_{h}\otimes T^{*}T(\xi_{h})) \\
& = & (\i\otimes \xi_{h}^{*})(\xi_{h}\otimes T^{*}T(\xi_{h})) \\
& = & \langle \xi_{h}, T^{*}T(\xi_{h})\rangle \xi_{h} \\
& = & \|T(\xi_{h})\|^{2}\xi_{h}
\end{eqnarray*}
and $\xi_{h}$ is an eigenvector for $A(T)^{*}A(T)$. If $A(T)^{*}A(T)$ is invertible, then
\begin{equation*}
(A(T)^{*}A(T))^{-1/2}\xi_{h} = \|T(\xi_{h})\|^{-1}\xi_{h}.
\end{equation*}
Thus, $A(T)\vert A(T)\vert^{-1}\xi_{h} = \xi_{h}\otimes \|T(\xi_{h})\|^{-1}T(\xi_{h})$ and 
\begin{eqnarray*}
\tau(M_{A(T)\vert A(T)\vert^{-1}} (x)) & = & \langle A(T)\vert A(T)\vert^{-1}\xi_{h}, (x\otimes 1)A(T)\vert A(T)\vert^{-1}\xi_{h}\rangle \\
& = & \langle \xi_{h}, x(\xi_{h})\rangle \\
& = & \tau(x).
\end{eqnarray*}
\end{proof}

Applying Lemma \ref{lem:technical1} to $M_{\delta} = A([P-Q]/2)$, and setting $\|x\|_{2} = \tau(x^{*}x)^{1/2}$ for $x\in D$, we can compute
\begin{equation*}
\|(m_{a} - M_{\gamma})(x)\|_{2}^{2} = \|M_{\delta}(x)\|_{2}^{2} = \tau(M_{\delta}(x)^{*}M_{\delta}(x))\leqslant \|\delta\|^{4}\|x\|_{2}^{2} \leqslant \epsilon^{4}\|x\|_{2}^{2}
\end{equation*}
i.e. $\|(m_{a} - M_{\gamma})(x)\|_{2}\leqslant \epsilon^{2}\|x\|_{2}$ and $M_{\gamma}$ approximates $m_{a}$ up to $\epsilon^{2}$ in $\mathcal{B}(L^{2}(D, \tau))$. We now only have to control $\|M_{\tilde{\gamma}} - M_{\gamma}\|_{\mathcal{B}(L^{2}(D, \tau))}$.

\begin{lemma}\label{lem:technical2}
$\tau((M_{\gamma}(x) - M_{\tilde{\gamma}}(x))^{*}(M_{\gamma}(x) - M_{\tilde{\gamma}}(x)))^{1/2}\leqslant 5\epsilon\tau(x^{*}x)^{1/2}$.
\end{lemma}

\begin{proof}
We have $\gamma = A((P+Q)/2)$. Thus, setting $T=((P+Q)/2)$ and observing that
\begin{equation*}
\left\|\left(T\xi_{h}- \frac{1}{\|T\xi_{h}\|}T\xi_{h}\right)\right\| = \left\|\xi_{h}\otimes\left(T\xi_{h}- \frac{1}{\|T\xi_{h}\|}T\xi_{h}\right)\right\| = \|(\gamma - \tilde{\gamma})\xi_{h}\|\leqslant \|\gamma - \tilde{\gamma}\| \leqslant \epsilon,
\end{equation*}
we can compute, again with Lemma \ref{lem:technical1},
\begin{eqnarray*}
\tau(M_{\gamma - \tilde{\gamma}}(x^{*}x)) & \leqslant & \langle (x^{*}x\otimes 1)(\gamma - \tilde{\gamma})\xi_{h}, (\gamma - \tilde{\gamma})\xi_{h}\rangle \\
& = & \left\langle (x^{*}x\otimes 1)\left(\xi_{h}\otimes \left(T\xi_{h}- \frac{1}{\|T\xi_{h}\|}T\xi_{h}\right)\right), \xi_{h}\otimes \left(T\xi_{h} - \frac{1}{\|T\xi_{h}\|}T\xi_{h}\right)\right\rangle \\
& = & \langle (x^{*}x)\xi_{h}, \xi_{h}\rangle \left\|T\xi_{h}- \frac{1}{\|T\xi_{h}\|}T\xi_{h}\right\|^{2} \\
& \leqslant & \epsilon^{2}\tau(x^{*}x) \\
& \leqslant & \epsilon^{2}\|x\|_{2}^{2}
\end{eqnarray*}
and $\|M_{\gamma - \tilde{\gamma}}(x)\|_{2}^{2} = \tau(M_{\gamma - \tilde{\gamma}}(x)^{*}M_{\gamma - \tilde{\gamma}}(x)) \leqslant \|\gamma - \tilde{\gamma}\|^{2}\tau(M_{\gamma - \tilde{\gamma}}(x^{*}x)) \leqslant \epsilon^{4}\|x\|^{2}_{2}$. Now, for all $x\in D$, we have
\begin{eqnarray*}
\|M_{\gamma}(x) - M_{\tilde{\gamma}}(x)\|_{2} & = & \|M_{\gamma}(x) - M_{\gamma + (\tilde{\gamma} -\gamma)}(x)\|_{2} \\
& \leqslant & \|M_{\tilde{\gamma} - \gamma}(x)\|_{2} + \tau((\tilde{\gamma} - \gamma)^{*}(x^{*}\otimes 1)\gamma\gamma^{*}(x\otimes 1)(\tilde{\gamma} - \gamma))^{1/2} \\
& & +\: \tau(\gamma^{*}(x^{*}\otimes 1)(\tilde{\gamma} - \gamma)(\tilde{\gamma} - \gamma)^{*}(x\otimes 1)\gamma)^{1/2} \\
& \leqslant & \epsilon^{2}\|x\|_{2} + (1+\epsilon)\tau((\tilde{\gamma} - \gamma)^{*}(x^{*}\otimes 1)(x\otimes 1)(\tilde{\gamma} - \gamma))^{1/2} \\
& & +\: \epsilon\tau(\gamma^{*}(x^{*}\otimes 1)(x\otimes 1)\gamma)^{1/2} \\
& \leqslant & \epsilon^{2}\|x\|_{2} + (1+\epsilon)\tau(M_{\tilde{\gamma} - \gamma}(x^{*}x))^{1/2} + \epsilon\tau(M_{\gamma}(x^{*}x))^{1/2} \\
& \leqslant & \epsilon^{2}\|x\|_{2} + (1+\epsilon)\epsilon\|x\|_{2} + \epsilon(1+\epsilon)\|x\|_{2} \\
& \leqslant & 5\epsilon\|x\|_{2}.
\end{eqnarray*}
Finally, $\|(m_{a} - M_{\tilde{\gamma}})(x)\|_{2}\leqslant 6\epsilon\|x\|_{2}$.
\end{proof}

We have to chek that this approximation also works in $\B(L^{2}(D, \tau)^{op})$, but this is straightforward.

\begin{lemma}\label{lem:technical3}
$M_{\tilde{\gamma}}$ also approximates $m_{a}$ up to $6\epsilon$ in $\mathcal{B}(L^{2}(D, \tau)^{op})$.
\end{lemma}

\begin{proof}
To estimate the opposite $L^{2}$-norm, one only needs to do all the previous computations exchanging $P$ and $Q$. Since they play symmetric rôles, we get the same result.
\end{proof}

Thus, we are able to approximate $m_{a}$ by a unital completely positive map in all the required norms. There remains only to check the $\tau$-invariance condition.

\begin{lemma}\label{lem:technical4}
The maps $M_{\tilde{\gamma}}$ and $m_{a}$ are $\tau$-preserving.
\end{lemma}

\begin{proof}
For $M_{\tilde{\gamma}}$, this comes from Lemma \ref{lem:technical1}. For $m_{a}$, this comes from the following computation
\begin{eqnarray*}
\tau(m_{a}(x)) & = & \langle \eta^{*}(x\otimes 1)\xi(\xi_{h}), \xi_{h}\rangle \\
& = & \langle (x\otimes 1)\xi(\xi_{h}), \eta(\xi_{h})\rangle \\
& = & \langle (x\otimes 1)A(P)(\xi_{h}), A(Q)(\xi_{h})\rangle \\
& = & \langle (x\otimes 1)(\xi_{h}\otimes P(\xi_{h}), \xi_{h}\otimes Q(\xi_{h})\rangle \\
& = & \langle x(\xi_{h}), \xi_{h}\rangle \langle P(\xi_{h}), Q(\xi_{h})\rangle \\
& = & \tau(x)\tau(m_{a}(1))
\end{eqnarray*}
and the fact that we have assumed $m_{a}$ to be unital.
\end{proof}

We are now ready to prove the theorem.

\begin{proof}[Proof of  Theorem \ref{thm:quantumfreeproduct}]
For each $i$, set $A_{i} = C_{\text{red}}(\mathbb{G}_{i})$ and $B_{i} = \mathcal{B}(L^{2}(\mathbb{G}_{i}))$. Consider a net $(a_{i, t})_{t}$ of finitely supported elements in $\ell^{\infty}(\h{\mathbb{G}}_{i})$ converging pointwise to the identity and such that $\limsup_{t} \|m_{a_{i, t}}\|_{cb} = 1$ and note that since $\h{\varepsilon}(a_{i, t})\rightarrow 1$ (because of the pointwise convergence assumption), we can, up to extracting a suitable subsequence, assume it to be non-zero and divide by it so that $m_{a_{i, t}}$ becomes unital. For any $0 < \epsilon < 1$, there is a $t(\epsilon)$ such that $\|m_{a_{i, t(\epsilon)}}\|_{cb}\leqslant 1+\epsilon$ (the same being automatically true for $m_{\h{S}(a_{i, t(\epsilon)})}$). Since
\begin{equation*}
\h{S}\circ\ast\circ\h{S}\circ\ast=\i,
\end{equation*}
we can replace $a_{i, t}$ by $(a_{i, t}+\h{S}(a_{i, t}))/2$ so that all the hypothesis of Lemma \ref{lem:cbapproximation} are satisfied. Then, by lemmata \ref{lem:cbapproximation}, \ref{lem:technical1}, \ref{lem:technical2}, \ref{lem:technical3} and \ref{lem:technical4} we get a unital completely positive approximation in completely bounded norm and in both $L^{2}$-norms. Applying Theorem \ref{thm:freecmap} proves that $\Lambda_{cb}(\ast_{i}A_{i}) = 1$. Since the original maps were all multipliers, the resulting finite-dimensional approximation is also implemented by multipliers and $\Lambda_{cb}(\ast_{i\in I}\h{\G}_{i}) = 1$.
\end{proof}

\begin{example}
Let $(G_{i})$ be any family of compact groups, then their duals in the sense of quantum groups are amenable (which implies that the Cowling-Haagerup constant is $1$ thanks to \cite[Thm 3.8]{tomatsu2006amenable}). Thus $\ast_{i}(C(G_{i}))$ is the dual of a non-cocommutative discrete quantum group with Cowling-Haagerup constant equal to $1$.
\end{example}

\begin{example}
The free orthogonal quantum groups $\widehat{A_{o}(F)}$ have Cowling-Haagerup constant equal to $1$ for any $F\in GL(2, \mathbb{C})$ such that $F\overline{F}=Id$ since they are amenable. Moreover, for any $F\in GL(2, \mathbb{C})$ the free unitary quantum group $\widehat{A_{u}(F)}$ is a quantum subgroup of $\mathbb{Z}\ast\widehat{A_{o}(F)}$, thus $\Lambda_{cb}(\widehat{A_{u}(F)})=1$. It follows that any reduced free product of some $2$-dimensional free quantum groups with duals of compact groups has Cowling-Haagerup constant equal to $1$.
\end{example}

\bibliographystyle{amsplain}
\bibliography{../../../quantum}

\providecommand{\bysame}{\leavevmode\hbox to3em{\hrulefill}\thinspace}
\providecommand{\MR}{\relax\ifhmode\unskip\space\fi MR }
\providecommand{\MRhref}[2]{%
  \href{http://www.ams.org/mathscinet-getitem?mr=#1}{#2}
}
\providecommand{\href}[2]{#2}
\begin{thebibliography}{10}

\bibitem{baaj1993unitaires}
S.~Baaj and G.~Skandalis, \emph{{Unitaires multiplicatifs et dualit{\'e} pour
  les produits crois{\'e}s de C*-alg{\`e}bres}}, Ann. Sci. \'Ecole Norm. Sup.
  \textbf{26} (1993), no.~4, 425--488.

\bibitem{bozejko1984herz}
Marek Bozejko and Gero Fendler, \emph{{Herz-Schur multipliers and completely
  bounded multipliers of the Fourier algebra of a locally compact group}},
  Boll. Un. Mat. Ital. A (6) \textbf{3} (1984), no.~2, 297--302.

\bibitem{daws2011multipliers}
M.~Daws, \emph{{Multipliers of locally compact quantum groups via Hilbert
  C*-modules}}, J. Lond. Math. Soc. \textbf{84} (2012), 385--407.

\bibitem{kraus1999approximation}
J.~Kraus and Z-J. Ruan, \emph{{Approximation properties for Kac algebras}},
  Indiana Univ. Math. J. \textbf{48} (1999), no.~2, 469--535.

\bibitem{ozawa2007weak}
N.~Ozawa, \emph{{Weak amenability of hyperbolic groups}}, Groups Geom. Dyn.
  \textbf{2} (2008), 271--280.

\bibitem{ozawa2010class}
N.~Ozawa and S.~Popa, \emph{{On a class of II$_{1}$ factors with at most one
  Cartan subalgebra}}, Ann. of Math. \textbf{172} (2010), no.~1, 713--749.

\bibitem{ozawa2010classbis}
\bysame, \emph{{On a class of II$_{1}$ factors with at most one Cartan
  subalgebra II}}, American J. Math. \textbf{132} (2010), no.~3, 841--866.

\bibitem{ricard2005khintchine}
E.~Ricard and Q.~Xu, \emph{Khintchine type inequalities for reduced free
  products and applications}, J. Reine Angew. Math. \textbf{599} (2006),
  27--59.

\bibitem{tomatsu2006amenable}
R.~Tomatsu, \emph{Amenable discrete quantum groups}, J. Math. Soc. Japan
  \textbf{58} (2006), no.~4, 949--964.

\bibitem{wang1995free}
S.~Wang, \emph{Free products of compact quantum groups}, Comm. Math. Phys.
  \textbf{167} (1995), no.~3, 671--692.

\bibitem{woronowicz1995compact}
S.L. Woronowicz, \emph{{Compact quantum groups}}, Sym{\'e}tries quantiques (Les
  Houches, 1995) (1998), 845--884.

\end{thebibliography}

\end{document}